\magnification 1250
\pretolerance=500 \tolerance=1000  \brokenpenalty=5000
\mathcode`A="7041 \mathcode`B="7042 \mathcode`C="7043
\mathcode`D="7044 \mathcode`E="7045 \mathcode`F="7046
\mathcode`G="7047 \mathcode`H="7048 \mathcode`I="7049
\mathcode`J="704A \mathcode`K="704B \mathcode`L="704C
\mathcode`M="704D \mathcode`N="704E \mathcode`O="704F
\mathcode`P="7050 \mathcode`Q="7051 \mathcode`R="7052
\mathcode`S="7053 \mathcode`T="7054 \mathcode`U="7055
\mathcode`V="7056 \mathcode`W="7057 \mathcode`X="7058
\mathcode`Y="7059 \mathcode`Z="705A
\def\spacedmath#1{\def\packedmath##1${\bgroup\mathsurround
=0pt##1\egroup$}
\mathsurround#1
\everymath={\packedmath}\everydisplay={\mathsurround=0pt}}
\def\nospacedmath{\mathsurround=0pt
\everymath={}\everydisplay={} } \spacedmath{2pt}

\font\eightrm=cmr8         \font\eighti=cmmi8
\font\eightsy=cmsy8        \font\eightbf=cmbx8
\font\eighttt=cmtt8        \font\eightit=cmti8
\font\eightsl=cmsl8        \font\sixrm=cmr6
\font\sixi=cmmi6           \font\sixsy=cmsy6
\font\sixbf=cmbx6\catcode`\@=11
\def\eightpoint{%
  \textfont0=\eightrm \scriptfont0=\sixrm \scriptscriptfont0=\fiverm
  \def\rm{\fam\z@\eightrm}%
  \textfont1=\eighti  \scriptfont1=\sixi  \scriptscriptfont1=\fivei
  \def\oldstyle{\fam\@ne\eighti}\let\old=\oldstyle
  \textfont2=\eightsy \scriptfont2=\sixsy \scriptscriptfont2=\fivesy
  \textfont\itfam=\eightit
  \def\it{\fam\itfam\eightit}%
  \textfont\slfam=\eightsl
  \def\sl{\fam\slfam\eightsl}%
  \textfont\bffam=\eightbf \scriptfont\bffam=\sixbf
  \scriptscriptfont\bffam=\fivebf
  \def\bf{\fam\bffam\eightbf}%
  \textfont\ttfam=\eighttt
  \def\tt{\fam\ttfam\eighttt}%
  \abovedisplayskip=9pt plus 3pt minus 9pt
  \belowdisplayskip=\abovedisplayskip
  \abovedisplayshortskip=0pt plus 3pt
  \belowdisplayshortskip=3pt plus 3pt
  \smallskipamount=2pt plus 1pt minus 1pt
  \medskipamount=4pt plus 2pt minus 1pt
  \bigskipamount=9pt plus 3pt minus 3pt
  \normalbaselineskip=9pt
  \setbox\strutbox=\hbox{\vrule height7pt depth2pt width0pt}%
  \normalbaselines\rm}\catcode`\@=12
\newcount\noteno
\noteno=0
\def\up#1{\raise 1ex\hbox{\sevenrm#1}}
\def\note#1{\global\advance\noteno by1
\footnote{\parindent0.4cm\up{\number\noteno}\
}{\vtop{\eightpoint\baselineskip12pt\hsize15.5truecm\noindent
#1}}\parindent 0cm}

\font\eightrm=cmr8     \font\sixrm=cmr6
\def\pc#1{\tenrm#1\sevenrm}
\def\tx{\kern-1.5pt -}
\def\cqfd{\kern 2truemm\unskip\penalty 500\vrule height 4pt
depth 0pt  width 4pt\medbreak} 
\def\no{n\up{o}\kern 2pt}
\def\rond{\kern 1pt{\scriptstyle\circ}\kern 1pt}
\frenchspacing
\font\tengoth=eufm10
  \font\sevengoth=eufm7
        \font\fivegoth=eufm5
\newfam\gothfam
\textfont\gothfam=\tengoth \scriptfont\gothfam=\sevengoth
\scriptscriptfont\gothfam=\fivegoth
\def\goth{\fam\gothfam\tengoth}%

\font\san=cmssdc10

\def\sym{\hbox{\san \char83}}
\def\ind{\par\hskip 1truecm\relax}
\def\indp{\par\hskip 0.5truecm\relax}
\def\Pic{\mathop{\rm Pic}\nolimits}
\def\pr{\mathop{\rm pr}\nolimits}
\def\mono{\lhook\joinrel\mathrel{\longrightarrow}}
\frenchspacing
\vsize = 25.3truecm
\hsize = 16.15truecm
\voffset = -.5truecm
\parindent=0cm
\baselineskip15pt
\overfullrule=0pt

\null{}\vskip1pt
\centerline{\bf  On the Chow ring of a K3 surface}
\smallskip
\smallskip \centerline{Arnaud {\pc BEAUVILLE} and Claire {\pc
VOISIN}}
\vskip1.2cm

{\bf Introduction}
\smallskip
\ind   An important algebraic invariant of a projective manifold $X$
is the  Chow ring $CH(X)$  of algebraic cycles on $X$ modulo linear
equivalence. It is  graded by the codimension  of cycles; the ring
structure comes from the  intersection product.
 For a surface we have
$$CH(X)={\bf Z}\oplus
\Pic(X)\oplus CH_0(X)\ ,$$where the group $CH_0(X)$
parametrizes $0$\tx cycles on $X$. While the structure of the
Picard group $\Pic(X)$ is well understood,  this is not the case
for $CH_0(X)$: if $X$ admits a nonzero holomorphic 2-form,  it is a
huge group, which cannot be parametrized by an algebraic variety
[M].
\ind Among the simplest examples of such surfaces are the K3
surfaces, which carry a nowhere vanishing holomorphic 2-form. In
this case $\Pic(X)$ is a lattice, while $CH_0(X)$ is very large; the
following result is therefore  somewhat  surprising: \smallskip {\bf
Theorem}$.-$ {\it Let $X$ be a {\rm K3} surface.
\indp{\rm a)} All points of $X$ which lie on some (possibly singular)
rational curve have the same class $c^{}_X$ in $CH_0(X)$.
\indp{\rm b)} The image of the intersection
product}$$\Pic(X)\otimes\Pic(X)\rightarrow CH_0(X)$${\it is
contained in ${\bf Z}\,c^{}_X.$
\indp{\rm c)} The second Chern class  $c_2(X)\in CH_0(X)$ is
equal to $24\,c^{}_X$}.
\smallskip
\ind The proof is elementary in the sense that it only appeals
 to simple
geometric constructions, based on the existence of sufficiently
many rational and elliptic curves on $X$. We prove a) and b) in
section 1; the proof of c), which is more involved, is given in
 section
2. If we represent the class $c^{}_X$ by a point $c$ of $X$, a key
property of this class is the formula
$$
    (x,x)-(x,c)-(c,x)+(c,c)=0\quad\hbox{\rm in}\quad
    CH_0(X\times X),$$
valid for any $x\in X$. In section 3 we discuss the importance of this formula and
its relation with property b) of the Theorem.
We  prove that an analogous formula holds when
$X$ is replaced by a hyperelliptic curve, but that it cannot hold for a generic
curve $C$ of genus $\ge 3$ -- we show that this would imply that $C$ is
algebraically equivalent to $-C$ in its Jacobian, contradicting a
result of Ceresa.
 \vskip1truecm {\bf 1. The image of the intersection
product}
\vglue3pt
\ind We work over the complex numbers. A {\it rational curve} on a
surface is  irreducible, but possibly singular.  If $V$ is an algebraic
variety and $p\in{\bf N}$, we denote by $CH_p(V)$ the group of $p$\tx
dimensional cycles on $V$ modulo rational equivalence; we put
$CH_p(V)_{\bf Q}=CH_p(V)\otimes{\bf Q}$.\smallskip (1.1) {\it  Proof  of} a)
{\it and} b) : Let
$R$ be a rational curve  on
$X$; it
 is the image of a generically injective map $j:{\bf
P}^1\rightarrow X$. Put $c^{}_R=j_*(p)$, where $p$ is an
arbitrary point of ${\bf P}^1$. For any divisor $D$
on $X$, we have in
$CH_0(X)$
$$R\cdot D=j_*j^*D=j_*(n\,p)= n\,c^{}_R\ ,\quad \hbox{with}\
n=\deg(R\cdot D)\ .$$
\ind Let $S$ be another
rational curve. If  $\deg(R\cdot S)\not=0$, the above
equality applied to $R\cdot S$ gives $c^{}_S=c^{}_R$ (recall that
$CH_0(X)$ is torsion free [R]). If  $\deg(R\cdot S)=0$,  choose an
ample divisor $H$; by a theorem of Bogomolov and Mumford
[M-M], $H$ is linearly equivalent to a sum of rational curves. Since
$H$ is connected, we can find a chain $R_0,\ldots,R_k$ of distinct
rational curves such that $R_0=R$, $R_k=S$ and $R_{i}\cap
R_{i+1}\not=\emptyset$ for $i=0,\ldots,k-1$. We conclude from
the preceding case that $c^{} _R=c^{} _{R_1}=\ldots=c^{}_S$.
Thus the class $c^{}_R$ does not depend on the choice of $R$:
this is assertion a) of the Theorem.
  Let us denote it by
$c^{}_X$.
\ind We have $R\cdot D =\deg(R\cdot D)\ c^{}_X$  for any divisor
$D$  and any rational curve $R$ on $X$. Since the group $\Pic(X)$
is spanned by the  classes of  rational curves (again by the
Bogomolov-Mumford theorem), assertion b) follows.\cqfd {\it
Remark} $1.2.-$ The  result (and the proof) hold  more generally
for  any  surface $X$ such that:
\indp a)  The Picard group of $X$ is
spanned by the  classes of  rational curves;
\indp b) There exists an ample divisor on $X$ which is a sum of
rational curves.
\ind This is the case when $X$ admits a non-trivial elliptic fibration
over ${\bf P}^1$ with a section, or for some particular surfaces like
Fermat surfaces in ${\bf P}^3$ with degree  prime to $6$ [S].
\smallskip  {\it
Remark} $1.3.-$  Let $A$ be an abelian surface. According to [Bl],
the image of the product map $\Pic(A)\otimes \Pic(A)\rightarrow
CH_0(A)$  has finite index, so the situation looks rather different
from  the K3 case. There is however an
 analogue to the Theorem. 
Let $\Pic^+(A)$ be the subspace of $\Pic(A)^{}_{\bf Q}$ fixed by
the action of the involution $a\mapsto -a$. We have a direct sum
decomposition $$\Pic(A)^{}_{\bf Q}=\Pic^+(A)\oplus \Pic^{\rm
o}(A)^{}_{\bf Q}\ ,$$so that
$\Pic^+(A)$ is canonically isomorphic to the image of $\Pic(A)^{}_{\bf Q}$ in
$H^2(A,{\bf Q})$. Now we claim that
{\it the image of the map $\mu :\Pic^+(A)\otimes
\Pic^+(A)\rightarrow CH_0(A)^{}_{\bf Q}$ is ${\bf Q}\,[0]$}. This
is a direct consequence of the decomposition of
$CH(A)^{}_{\bf Q}$ described in [B]: let $k$ be an integer $\ge 2$,
and let ${\bf k}$  be the multiplication by $k$ in $A$.
We have ${\bf k}^*D=k^2D$ for any element $D$ of  $\Pic^+(A)$,
thus ${\bf k}^*c^{} =k^4c^{} $ for any element $c^{} $ in the image
of  $\mu $; but the latter property characterizes the multiples of
$[0]$.\cqfd
\medskip\ind (1.4)
The cycle class
$c^{}_X$ has some remarkable properties that we will  investigate in the next
section.
 Let us observe first that for any irreducible curve $C$ on $X$,
 there is a  rational curve $R\not= C$ which intersects $C$; thus we
can represent $c^{}_X$ by the class of a point $c\in C$ (namely any
point of $C\cap R$).\smallskip
\ind We will need a more subtle property of $c^{}_X$. Let us first prove
a lemma:
\smallskip {\bf Lemma 1.5}$.-$ {\it Let $E$ be an elliptic curve,
$x,y$ two points of $E$. Then} $$(x,x)-(x,y)-(y,x)+(y,y)=0\qquad
\hbox{\it in}\qquad CH_0(E\times E)\ .$$
\ind Since the divisors $[x]-[y]$ generate the group
$\Pic^{\rm o}(E)$, this is equivalent to the formula
$\pr_1^*D\cdot
\pr_2^*D=0$ in $CH_0(E\times E)$ for every $D$ in $\Pic^{\rm
o}(E)$.\smallskip {\it Proof} :  Put $\xi=(x,x)-(x,y)-(y,x)+(y,y)$. Then
$2\xi$ is the pull-back of a $0$\tx cycle $\eta=(x,x)+(y,y)-2(x,y)$
on the second symmetric product $\sym^2E$. The addition map
$a:\sym^2E\rightarrow E$ is a ${\bf P}^1$\tx fibration; this implies
that the push-down map
$a_*:CH_0(E\times E)\rightarrow CH_0(E)$ is an isomorphism. Since
$a_*\eta=0$, we have $\eta=0$, hence $2\xi=0$. On the other
hand $\xi$ has degree $0$ and its image in the Albanese variety of
$E\times E$ is zero, so $\xi=0$ by Rojtman's result.\cqfd

\smallskip  {\bf Proposition 1.6}$.-$ {\it Let
$\Delta$  be the diagonal embedding of
$X$ into $X\times X$.

\indp {\rm a)} For every $\alpha\in CH_1(X)$, we have
$$ \Delta_*\alpha=\pr_1^*\alpha\cdot
\pr_2^*c^{}_X+\pr_1^*c^{}_X\cdot
\pr_2^*\alpha\quad \hbox{in}\quad CH_1(X\times X)\ ;$$
\indp {\rm b)} For every $\xi\in CH_0(X)$, we have
$$\Delta_*\xi=\pr_1^*\xi\cdot\pr_2^*c^{}_X+\pr_1^*c^{}_X\cdot
\pr_2^*\xi-(\deg \xi)\,\Delta_*c^{}_X\quad
\hbox{in}\quad CH_0(X\times X)\ .$$}\par
{\it Proof} : a) Since both sides are additive in $\alpha$, it is
enough to check this relation when $\alpha$ is the class of a rational
curve; in that case it follows from the fact that the diagonal of ${\bf
P}^1\times {\bf P}^1$ is linearly equivalent to
${\bf P}^1\times \{0\}+\{0\}\times {\bf P}^1$.
\ind b) Again both sides are additive in $\xi$, so we may assume
that $\xi$ is the class of a point $x\in X$. The Bogomolov-Mumford
theorem tells us that $x$ lies on the image of a curve $E$ of genus
$\le 1$; by (1.4) we can represent $c^{}_X$ by a point $c\in E$.
We have $(x,x)-(x,c)-(c,x)+(c,c)=0$ in $CH_0(E\times E)$ by
lemma 1.5 (the case when $E$ is rational is trivial); by push-down this
gives the same formula in $CH_0(X\times X)$.\cqfd
\vskip1truecm
\font\grit=cmmib10%\font\goth=eufm7
\def\cg{\hbox{\grit \char 99}}
\def\gox{{\scriptstyle\goth X}}
{\bf 2. The formula ${\bf \cg_2(X)=24\,\cg_X^{}}$}\vglue3pt
\ind(2.1) Let $c$ be a point of $X$ lying on some rational curve. We
will denote  by
$(x,x,x)$, $(x,x,c)$, $(x,c,c)$, etc...
the classes in $CH_2(X\times X\times X)$ of the image of $X$ by the
maps
$\  x\mapsto (x,x,x)\ ,\  x\mapsto(x,x,c)\ , \
x\mapsto(x,c,c)\ ,\ {\rm etc...}$
With this notation we
have the following key result:
\smallskip
{\bf Proposition 2.2}$.-$ {\it The
cycle}$${\gox}=(x,x,x) - (c,x,x)- (x,c,x)- (x,x,c) +
(x,c,c) + (c,x,c)+ (c,c,x)$$ {\it is zero in
$CH_2(X\times X\times X)^{}_{\bf Q}$}.
\medskip
{\bf Corollary 2.3}$.-$ {\it Let
$\Delta$, $i_c$ and $j_c$ be the
maps of $X$ into $X\times
X$ defined by $\Delta(x)=(x,x)$, $i_c(x)= (x,c)$ and
$j_c(x)=(c,x)$.  For every $\xi$ in $CH_2(X\times
X)$, we have}
$$\Delta^*\xi=i_c^*\xi+j_c^*\xi+n\,c,\qquad\hbox{\it
with}
\qquad n=\deg(\Delta^*\xi-i_c^*\xi-j_c^*\xi)\ .$$
\smallskip
\ind From this formula we recover part b) of the
Theorem by taking $\xi={\rm pr}_1^*\alpha\cdot {\rm
pr}_2^*\beta$, with $\alpha,\beta\in\Pic(X)$, and
we get part c) by taking for $\xi$ the class of the
diagonal, so that
$\Delta^*\xi=c_2(X)$. \medskip
(2.4) {\it Proof of the Corollary} : We will denote by
$p_i$, for
$1\le i\le 3$, the projection of $X\times X\times X$ onto
the $i$\tx th factor, and by $p_{ij}$, for
$1\le i<j\le 3$, the projection
$(x_1,x_2,x_3)\mapsto (x_i,x_j)$.
  \ind Let us compute $p_3{}_*(\gox\cdot
p_{12}^*\xi)$. Let $\delta:X\rightarrow X\times X\times
X$  be the  map $x\mapsto(x,x,x)$. We have $p_3\rond
\delta ={\rm Id}_X$ and $p_{12}\rond\delta =\Delta$,
hence
$$p_3{}_*((x,x,x)\cdot
p_{12}^*\xi)=p_3{}_*\delta_*(\delta
^*p_{12}^*\xi)=\Delta^*\xi\ .$$
\ind The same argument applied to the maps
$x\mapsto (c,x,x)$, $x\mapsto (x,c,x)$, ... ,
gives
$$\nospacedmath\displaylines{p_3{}_*((c,x,x)\cdot
p_{12}^*\xi) =i_c^*\xi \qquad
p_3{}_*((x,c,x)\cdot p_{12}^*\xi)=j_c^*\xi\cr
p_3{}_*((x,x,c)\cdot
p_{12}^*\xi)=\deg(\Delta^*\xi)\cdot c \qquad
p_3{}_*((x,c,c)\cdot p_{12}^*\xi)=\deg(i_c^*\xi)\cdot
c\cr
p_3{}_*((c,x,c)\cdot p_{12}^*\xi)=\deg(j_c^*\xi)\cdot c
\qquad p_3{}_*((c,c,x)\cdot p_{12}^*\xi)=0\ ,}$$hence
our formula.\cqfd
{\it Remark} $2.5.-$ One also recovers Proposition 1.6 b)
by restricting for each $x\in X$ the class $\gox$ to
$X\times X\times \{x\}\subset X\times X\times X$.
\medskip
\ind For the proof of the proposition we will need two
results on products of elliptic curves. Let $F$ be an
elliptic curve over an arbitrary field. We denote by
 $\Pic(F^3)^{\rm inv}$ the subgroup of elements of $\Pic(F^3)_{\bf Q}$ which
are invariant under permutations of the factors and under the involution
$(-1_{F^3})$. We keep the notation of  (2.1).
\smallskip
{\bf Lemma 2.6}$.-$ a) {\it
 The cycle class
$${\goth
v}=(u,u,u)-(0,u,u)-(u,0,u)-(u,u,0)+(u,0,0)+(0,u,0)+(0,0,u)$$
in $CH_1(F^3)^{}_{\bf Q}$ is zero}.
 \indp {\rm b)} {\it The
 divisors
$\ \alpha^{}_F=\displaystyle \sum_ip_i^*0\ $ and $\ \beta^{}_F
=\displaystyle \sum_{i<j}p_{ij}^*\Delta\
$ form a basis  of}
$\Pic(F^3)^{\rm inv}$.

\smallskip
{\it Proof} : a) The class ${\goth v}$ is symmetric,
hence comes from a cycle class $\bar {\goth v}$ in the
third symmetric product $\sym^3 F$. This variety  is a
${\bf P}^2$\tx bundle over
$F$, through the addition map $a:\sym^3
F\rightarrow F$. Thus we have  $CH_1(\sym^3
F)=a^*\Pic(F)\cdot h\oplus {\bf Z}h^2$, where $h$ is any
divisor class on $\sym^3 F$ which induces on a fibre
$a^{-1} (u)\cong {\bf P}^2$ the class of a line.
\ind Write $\bar {\goth v}=(a^*d)\cdot h+nh^2$. We have
$n=\deg(\bar {\goth v}\cdot a^*0)=3^2-3\cdot 2^2+3\cdot
1=0$, hence
$d=a_*(\bar {\goth v}\cdot h)$. We can represent $h$ by
the image of the divisor $p_1^*0$ in $F\times F\times
F$; since ${\goth v}\cdot p_1^*0=0$, we get $d=0$ and
finally
${\goth v}=0$.\smallskip
\ind b) As above we have $\Pic(\sym^3 F)=a^*\Pic(F)
\oplus {\bf Z}h$. Taking
the invariants under $(-1_{F^3}) $ we see that
$\Pic(F^3)^{\rm inv}$ has rank $2$. Thus it
suffices to prove that the divisors $\alpha^{}_F$
and $\beta^{}_F$ are not proportional in
$\Pic(F^3)$; but their restriction to $F^2$ (embedded in $F^3$
by $(u,v)\mapsto (u,v,0)$) are 
$p_1^*0+p_2^*0$ and $\Delta+p_1^*0+p_2^*0$, which are clearly
non-proportional.\cqfd
\medskip

\ind (2.7) We now begin the proof of the
Proposition. It will make our life easier to assume that
$\Pic(X)$ is generated by an ample divisor class $H$;
the general case will follow by specialization (see
[SGA6],  X.7.14). By the Bogomolov-Mumford theorem, we can
find in the linear system $|H|$ a one-dimensional family
$(E'_b)^{}_{b\in B}$ of (singular) elliptic curves; that
is, we can find a  surface $E$ with a fibration
$p:E\rightarrow B$ onto a smooth curve, with general
fibre a smooth curve $E_b$ of genus $1$, 
and a generically finite map
$\pi :E\rightarrow X$ which maps each fibre $E_b$ of $p$
birationally onto the singular curve $E'_b$. Passing to a
covering of
$B$ if necessary, we may assume that:
\indp a)
$p$ has a section $0:B\rightarrow E$;
\indp b) The curve $\pi (0_B)$ is {\it rational}.\par
(To see b), replace
$B$ by a component of $\pi ^{-1} (R)$, where $R$ is a
rational curve on $X$ not contained in any $E'_b$.)
\ind Note that because of the assumption on $\Pic(X)$
every fibre $E_b$ is irreducible.
\bigskip
\ind (2.8) Using again the  notation of (2.1), we
consider  on the fibre product
$E^3_B=$ $E\times _BE\times _BE$ the cycle class
$$\nospacedmath\displaylines{\quad {\goth u}=(u,u,u) -
(0_{pu},u,u)- (u,0_{pu},u)- (u,u,0_{pu}) +\hfill\cr\hfill
(u,0_{pu},0_{pu}) + (0_{pu},u,0_{pu})+
(0_{pu},0_{pu},u)\ .\quad }$$
\ind For $b\in B$, the class in $CH_2(X\times X\times
X)$ of the cycle
$$\{c\}\times E'_b\times E'_b\ +\ E'_b \times \{c\} \times E'_b\
+\ \{c\}\times E'_b\times E'_b $$
does not depend on $b$, since the curves $E'_b$ all
belong to the same linear system $|H|$; let
us denote it by ${\goth z}$.
Let $\pi ^3:E^3_B\rightarrow X^3$ be the morphism
deduced from $\pi $.
\medskip
{\bf Lemma 2.9}$.-$
 {\it The class $\pi ^3_*({\goth u})$ is proportional to
${\goth z}$}.
\smallskip {\it Proof} : By lemma 2.6.a), the
restriction of ${\goth u}$ to the generic fibre of the fibration
$E^3_B\rightarrow B$ is zero.
It follows that ${\goth u}$ is a sum of cycles of the form
$i_b{}_*D_b$, where $i_b$ is the inclusion of $E_b^3$ into
$E^3_B$ and $D_b$ a (Weil) divisor on $E_b^3$ [Bl-S].
\ind The involution $\sigma $ of $E$ which coincides with
$u\mapsto -u$ on each smooth fibre gives rise to an
involution $\sigma^3 $ of $E^3_B$ which commutes with the action
 of ${\goth S}_3$ by permutations  of the factors. The cycle
${\goth u}$ is invariant by this action of ${\goth S}_3\times {\bf
Z}/2$. By averaging on this group we may choose the above divisor
classes
$D_b$ in the invariant subgroup of $CH_2(E^3_b)^{}_{\bf Q}$.
We want to prove that each cycle class $i_b{}_*D_b$ is pushed
down to a multiple of ${\goth z}$ by
$\pi ^3$.
\ind Assume first that the curve $E_b$ is smooth. By
lemma 2.6.b)  the class $D_b$ is a ${\bf Q}$\tx linear
combination of
$\alpha^{}_{E_b}$ and $\beta^{}_{E_b}$. By (2.7.b) $\pi(0_b)$ is
linearly equivalent to $c$, thus we have $\pi
^3_*(i_b{}_*\,\alpha^{}_{E_b})={\goth z}$. The cycle $\pi
^3_*(i_b{}_*\,\beta^{}_{E_b})$ is the sum of $(\Delta_*E'_b)\times
E'_b$ and of the two cycles obtained by permutation of the
factors. Now using lemma 1.4  this class is equivalent
to $2{\goth z}$, hence the result in this case.
\ind If $E_b$ is singular, its normalization $\widetilde{E}_b$
is a smooth rational curve, and we have a surjective
homomorphism
$CH_2(\widetilde{E}_b^3)^{}_{\bf Q}\rightarrow
CH_2(E_b^3)^{}_{\bf Q}$. The ${\goth S}_3$\tx invariant part of
$CH_2(\widetilde{E}_b^3)^{}_{\bf Q}$ is spanned by the divisor
$\alpha^{}_{\widetilde{E}_b}=\sum p_i^*0$
 which again maps to a cycle linearly equivalent to
${\goth z}$ under $\pi ^3$.\cqfd

 {\bf Lemma 2.10}$.-$ {\it Let $d=\deg \pi $, and let $\gox$ be
the cycle  class  defined in} (2.2). {\it Then} $$\pi
^3_*({\goth u})=d\,\gox\ .$$ \smallskip {\it Proof} : We
compute the images under $\pi ^3_*$ of the cycles which appear
in the definition of $\gox$.
\ind a) We have $\pi ^3_*(u,u,u)=d\,(x,x,x)$ in $CH_2(X\times
X\times X)$.
\ind b) Let $\Gamma \i X\times X$ be the image of the
surface $(u,0_{pu})$ (that is, the graph of $0\rond p$) in
$E\times E$. We have
$\pi ^3_*(u,u,0_{pu})=p_{12}^*\Delta\cdot p_{23}^*\Gamma $.
The normalization $\widetilde{R}$ of $R=\pi
(0_B)$ is a smooth rational curve (2.7.b). Since our
cycle $\Gamma
$ is supported by $X\times R$, it comes from a divisor $\Gamma
_0$ in
$X\times \widetilde{R}$. Such a divisor is of the form $D\times
\widetilde{R}+ mX\times \{r\}$ for some divisor $D$ on $X$,
some point $r\in\widetilde{R}$ and some integer $m$; this
integer is equal to the degree of
$\Gamma _0$ over $X$, that is $d$. Therefore $\Gamma $
is linearly equivalent to $d\,(X\times c)+ D\times R$; since we
assume $\Pic(X)={\bf Z}$ we have $D\times R=a\,E'_b\times E'_b$
in $CH_2(X\times X)$ for some integer $a$ and any $b\in B$.
Intersecting with
$p_{12}^*\Delta$ we get
$$\pi ^3_*(u,u,0_{pu})=d\,(x,x,c)+a\,(\Delta_*E'_b)\times
E'_b\ .$$
\ind  c) We have $\pi
^3_*(u,0_{pu},0_{pu})=p_{12}^*\Gamma \cdot p_{23}^*\Delta $;
reasoning as in b) we  find
$$\pi ^3_*(u,0_{pu},0_{pu})=d\,(x,c,c)+a\,E'_b\times
(\Delta_*E'_b)\ .$$
\ind d) The lemma follows by  permuting and summing.\cqfd
\smallskip
\ind (2.11) Therefore $\gox=d^{-1} \,\pi ^3_*({\goth u})$ is
proportional to the effective cycle ${\goth z}$ (Lemma 2.9). On the
other hand
$\gox$ is homologically trivial: this
follows from the K\"unneth formula and the fact that the
cycles $p_{ij}{}_*\gox$ are identically zero. Thus we obtain
$\gox=0$, which concludes the proof of the Proposition.\cqfd
\vskip1truecm
{\bf 3. ${\bf 0}$\tx cycles on a product}\smallskip
\ind (3.1) The  cycle $c_X^{}$ has two remarkable properties, 
namely the intersection property b) of the Theorem, and the diagonal
property $(x,x)-(x,c)-(c,x)+(c,c)=0$ in
$CH_0(X\times X)$ for any  $x\in X$. These two properties may seem unrelated.
However we can rephrase them in the following way, which shows that they are in
some sense dual to each other: since the Picard group of 
$X$ is isomorphic to its N\'eron-Severi group,
the degree $1$ zero-cycle $c_X^{}$
provides a splitting of
$CH(X)$ as
 $$CH(X)=CH(X)_{\rm hom}\oplus H,$$
 where $CH(X)_{\rm hom}$ is the subgroup of $0$\tx cycles homologous to
$0$, and $H$  the image of $CH(X)$ into
$H^*(X)$ via the
 cycle map. This splitting induces the splitting
 $$\nospacedmath\displaylines{CH(X)\otimes
CH(X)=\hfill\cr \hfill (CH(X)_{\rm hom}\otimes CH(X)_{\rm hom})\,
 \oplus\, (CH(X)_{\rm hom}\otimes H)\, \oplus\,  (H\otimes CH(X)_{\rm hom})
\,\oplus\,  (H\otimes H)}$$
 of $CH(X)\otimes CH(X)$. It is immediate to see that it
 induces one
 on the image $CH(X\times X)_{\rm dec}$ of
 $CH(X)\otimes CH(X)$ in $CH(X\times X)$.
 We can see these decompositions as giving  gradings on $CH(X)$ and 
 $CH(X\times X)_{\rm dec}$. (Here it is natural from the point of view
 of the Bloch-Beilinson conjectures to assign the degree
 $0$ to $H$ and the degree $2$ to $CH(X)_{\rm hom}$ since 
 our surface is regular.)
 Then the intersection property b) says that
 if $\Delta:X\rightarrow X\times X$ is the diagonal embedding,
 the homomorphism
 $$\Delta^*:CH(X\times X)_{\rm dec}\rightarrow CH(X)$$
 is compatible with the gradings,
 while the diagonal relations a) and b) of  Proposition 1.6 say that for $p<2$
 the homomorphism
 $$\Delta_*:CH_p(X)\rightarrow CH_p(X\times X)$$
 takes values in 
 $CH_p(X\times X)_{\rm dec}$ and
 is also compatible with the gradings.

\medskip
\ind We are now going to  investigate the corresponding diagonal
property for a curve.
\smallskip
{\bf Proposition 3.2}$.-$ {\it Let $C$ be a hyperelliptic curve, and
$w$ a Weierstrass point of $C$.  For any $x\in C$, we have
$$(x,x)-(x,w) -(w,x)+(w,w)=0 \quad\hbox{\it in}\quad
CH_0(C\times C)\ .$$}
\hskip1truecm (Note that the class of $w$ is well-defined in
$CH_0(C)_{\bf Q}$.)
\smallskip {\it Proof} : Let $J$ be the Jacobian variety of $C$;
choose an  Abel-Jacobi embedding $C\mono J$. The induced map
$C\times C\rightarrow J\times J$ is an Albanese map for
$C\times C$.
\ind The subgroup  of
 degree $0$ cycles  in $CH_0(C\times C)$ maps onto the Albanese
variety
$J\times J$; let $T(C\times C)$ be the kernel of this map.
 The surjective  map
$$CH_0(C)\otimes CH_0(C)\longrightarrow CH_0(C\times C)$$
 induces a surjective map
$$J\otimes J\longrightarrow
T(C\times C)\ .$$ \ind  Let $\iota$ be the hyperelliptic
involution of $C$; since $\iota$ acts as $(-1)$ on $J$, we see that
{\it the involution $(\iota,\iota)$ of $C\times C$ acts
trivially on} $T(C\times C)$.
\ind Let ${\goth c}:J\rightarrow CH_0(C\times
C)$ be the homomorphism defined by $${\goth
c(\alpha)}=\Delta_*\alpha-\pr_1^*\alpha\cdot
\pr_2^*w-\pr_1^*w\cdot \pr_2^*\alpha\ .$$
 The
cycle ${\goth c}(\alpha)$ is of degree zero, and its image in
$J\times J$  is $(\alpha,\alpha)-(\alpha,0)-(0,\alpha)$ $=0$, hence
it is  invariant under
$(\iota,\iota)$. On the other hand we have $$(\iota,\iota)^*{\goth
c}(\alpha)={\goth c}(\iota^*\alpha)={\goth c(-\alpha)}=-{\goth
c(\alpha)}\ .$$ Therefore $2{\goth c}(\alpha)=0$, and actually
${\goth c}(\alpha)=0$ by Rojtman's result. Applying this to
$\alpha=[x]-[w]$ gives the result.\cqfd
\ind In contrast, we have now :
\smallskip 
{\bf Proposition 3.2}$.-$ {\it Let $C$ be a general curve of genus
$\ge 3$. There  exists no divisor $c$ on $C$ such that the
$0$\tx cycle\note{Here $(x,c)$ stands for the $0$\tx cycle
$\pr_1^*x\cdot\pr_2^*c$}
$(x,x)-(x,c)-(c,x)$ in
$CH_0(C\times C)$ is independent of $x\in C$}.

\smallskip {\it Proof} : As above the hypothesis on $c$ is
equivalent to the relation
$$\Delta_*\alpha=\pr_1^*\alpha\cdot \pr_2^*c + \pr_1^*c\cdot
\pr_2^*\alpha$$ for all $\alpha$ in $J$.
Applying $\pr_1{}_*$ we observe that this formula implies $\deg
c=1$.
\ind Put $c'= (x,c)+(c,x)-(x,x)$, and assume that this class in
$CH_0(C\times C)$ is independent of $x$. With the notation of
(2.1), we consider in
$CH_1(C\times C\times C)_{\bf Q}$  the cycle
$${\goth z}=(x,x,x)-(c,x,x)-(x,c,x)-(x,x,c)+(c,c,x)+(c,x,c)+(x,c')\ .
$$
Our hypothesis ensures that the
restriction of
${\goth z}$ to the generic fibre of
$p_1$ is zero. As in [Bl-S] we conclude that ${\goth z}$ is a sum of
1-cycles of the form $i_b{}_*D_b$, where $i_b:C\times
C\rightarrow C\times C\times C$ is the embedding $(x,y)\mapsto
(b,x,y)$ and $D_b$ is a divisor on $C\times C$.
\ind Let us now work in the group $A_1(C\times C\times C)_{\bf
Q}$ of cycles modulo algebraic equivalence. In this group
 the class of
$i_b{}_*D$, for $D\in CH_1(C\times C)_{\bf Q}$, is independent of
$b\in C$; thus we can write
${\goth z}=i_b{}_*D$ for some fixed $b\in C$ and some divisor $D$
in $C\times C$. Since $p_{12}\rond i_b={\rm Id}_{C\times
C}$ we have $D=p_{12}{}_*{\goth z}$.
\ind Now  the cycle ${\goth z}$ is homologically
trivial: as in (2.11) it suffices to check this for the projections
$p_{ij}{}_*{\goth z}$ on $C\times C$, and this is
straightforward. Thus, the divisor
$D$ is homologically, and therefore algebraically, trivial
in $C\times C$;
 we conclude
that
${\goth z}$ is zero in $A_1(C\times C\times C)_{\bf Q}$.
\ind Now let $J$ be the Jacobian variety of $C$, and
$\alpha:C\rightarrow J$ the Abel-Jacobi map which maps a
point $x$ of $C$ to the divisor class  $[x]-c$; we will identify $C$
with its image under $\alpha$. Let $\alpha^3:C^3\rightarrow J$ be
the map deduced from
$\alpha$. We have
$$(\alpha^3)_*({\goth z})={\bf 3}_*C\,-\,3({\bf
2}_*C)\,+\,3C=0\quad{\rm in}\quad A_1(J)_{\bf Q}\ ,$$where
${\bf k}$
 denotes the multiplication   in $J$ by the integer $k$.
\ind According to [B] we have a decomposition
$$A_1(J)_{\bf Q}=A_1(J)_0\oplus \cdots \oplus A_1(J)_{g-1}\
,$$where ${\bf k}_*$ acts by multiplication by $k^{2+s}$ on
$A_1(J)_s$. Since $3^\ell -3\cdot 2^\ell +3>0$ for $\ell \ge 3$,
the above equality  implies that  the components of the
$1$\tx cycle
$C$ in $A_1(J)_i$ are zero for $i\ge 1$, that is,  $[C]\in A_1(J)_0$.
Taking $k=-1$ we see that $C$ is algebraically equivalent to $-C$;
this contradicts the result of Ceresa [C].\cqfd

\vskip2cm
\centerline{ REFERENCES} \vglue15pt\baselineskip12.8pt
\def\num#1{\smallskip\item{\hbox to\parindent{\enskip
[#1]\hfill}}}
\parindent=1.3cm
\num{B} A. {\pc BEAUVILLE}: {\sl 	Sur l'anneau de Chow d'une
vari\'et\'e ab\'elienne.} Math. Annalen {\bf 273} (1986),
647--651.
\num{Bl} S. {\pc BLOCH}: {\sl Some elementary theorems about
algebraic cycles on Abelian varieties}. Invent. Math. {\bf 37}
(1976),  215--228.
\num{Bl-S} S. {\pc BLOCH}, V. {\pc SRINIVAS}:
{\sl Remarks on correspondences and algebraic cycles}. Amer. J. Math.
 {\bf 105} (1983), 1235--1253.
\num{C} G. {\pc CERESA}: {\sl $C$ is not algebraically equivalent to
$C\sp{-}$ in its Jacobian}. Ann. of Math. {\bf 117} (1983),
285--291.
%\num{G-G} M. {\pc GREEN}, Ph. {\pc GRIFFITHS}:
%{\sl An interesting $0$\tx cycle}, preprint 2001.
\num{M} D. {\pc MUMFORD}: {\sl Rational equivalence of
$0$\tx cycles on surfaces}.  J. Math. Kyoto Univ. {\bf 9} (1968),
195--204.
\num{M-M} S. {\pc MORI}, S. {\pc MUKAI}: {\sl Mumford's theorem
on curves on $K3$ surfaces}. Algebraic Geometry (Tokyo/Kyoto
1982), LN {\bf 1016},  351--352; Springer-Verlag (1983).
\num{R} A. A. {\pc ROJTMAN}: {\sl The torsion of the group of
$0$-cycles  modulo rational equi\-valence}.  Ann. of Math. {\bf
111} (1980), 553--569.
\num{S} T. {\pc SHIODA}: {\sl On the Picard number of a Fermat
surface}. J. Fac. Sci. Univ. Tokyo {\bf 28} (1982), 725--734.
\num{SGA6} {\sl Th\'eorie des intersections et th\'eor\`eme de
Riemann-Roch}.  S\'eminaire de G\'eom\'etrie Alg\'ebrique
du Bois-Marie 1966--1967 (SGA 6). Dirig\'e par
P.~Ber\-thelot, A. Grothendieck et L. Illusie.  Lecture
Notes in Math. {\bf  225}, Springer-Verlag, Berlin-New
York (1971).
\vskip1cm
\def\pc#1{\eightrm#1\sixrm}

$$\kern-0.1cm\vtop{\eightrm\baselineskip12pt\hbox to 5cm{\hfill Arnaud
{\pc BEAUVILLE}\hfill}\smallskip
\hbox to 5cm{\hfill Institut Universitaire de France\hfill}\vskip-2pt
\hbox to 5cm{\hfill \&\hfill}\vskip-2pt
 \hbox to 5cm{\hfill Laboratoire J.-A. Dieudonn\'e\hfill}
 \hbox to 5cm{\sixrm\hfill UMR 6621 du CNRS\hfill}
\hbox to 5cm{\hfill {\pc UNIVERSIT\'E DE}  {\pc NICE}\hfill}
\hbox to 5cm{\hfill  Parc Valrose\hfill}
\hbox to 5cm{\hfill F-06108 {\pc NICE} Cedex 2\hfill}}
\kern2,9cm
\vtop{\eightrm\hbox to 5cm{\hfill Claire {\pc
VOISIN}\hfill}\smallskip
 \hbox to 5cm{\hfill Institut de Math\'ematiques de
          Jussieu\hfill} \hbox to 5cm{\sixrm\hfill UMR 7586 du
CNRS\hfill}
\hbox to 5cm{\hfill   Case 7012 \hfill}
\hbox to 5cm{\hfill  2 place Jussieu\hfill}
\hbox to 5cm{\hfill F-75251 {\pc PARIS} Cedex 05\hfill}}$$
\end